\magnification=\magstep1
\input amstex
\documentstyle{amsppt}

\font\tencyr=wncyr10 
\font\sevencyr=wncyr7 
\font\fivecyr=wncyr5 
\newfam\cyrfam \textfont\cyrfam=\tencyr \scriptfont\cyrfam=\sevencyr
\scriptscriptfont\cyrfam=\fivecyr
\define\hexfam#1{\ifcase\number#1
  0\or 1\or 2\or 3\or 4\or 5\or 6\or 7 \or
  8\or 9\or A\or B\or C\or D\or E\or F\fi}
\mathchardef\Sha="0\hexfam\cyrfam 58

\define\Primes{\frak{Primes}}

\define\defeq{\overset{\text{def}}\to=}
\define\ab{\operatorname{ab}}
\define\pr{\operatorname{pr}}
\define\Gal{\operatorname{Gal}}
\def \isom {\overset \sim \to \rightarrow}

\define\Spec{\operatorname{Spec}}
\define\id{\operatorname{id}}

\define\Ker{\operatorname{Ker}}

\def \Sec{\operatorname {Sec}}

\def\genus{\operatorname{genus}}
\def\sol{\operatorname{sol}}
\def \Res{\operatorname{Res}}
\def \all{\operatorname{all}}

\def \tor{\operatorname{tor}}

\def \Hom{\operatorname{Hom}}
\def \and{\operatorname{and}}

\def \nd{\operatorname{nd}}

\NoRunningHeads
\NoBlackBoxes
\topmatter

\title
A local-global principle for torsors under geometric prosolvable fundamental groups
\endtitle

\author
Mohamed Sa\"\i di
\endauthor

\abstract We prove a {\it local-global principle for torsors under the prosolvable geometric fundamental group of
a hyperbolic curve over a number field}. 
\endabstract

\endtopmatter

\document

\subhead
\S 0. Introduction
\endsubhead
Let $k$ be a field of characteristic $0$, and $X\to \Spec k$ a proper, smooth, and geometrically connected 
{\it hyperbolic curve} over $k$ (i.e., $\genus X\ge 2$). Let $\eta$ be a geometric point of $X$ with value in the generic point of $X$.
Thus, $\eta$ determines naturally an algebraic closure $\bar k$ of $k$, and a geometric point $\overline \eta$
of $\overline X\defeq X\times _k\bar k$.
There exists a canonical exact sequence of profinite groups (cf. [Grothendieck], Expos\'e IX, Th\'eor\`eme 6.1)
$$1\to \pi_1(\overline {X},\overline \eta)\to \pi_1(X, \eta) @>>> G_k\to 1.$$
Here, $\pi_1(X, \eta)$ denotes the {\it arithmetic \'etale fundamental group} of $X$ with base
point $\eta$, $\pi_1(\overline {X},\overline \eta)$ the {\it \'etale fundamental group} of $\overline {X}$ with base
point $\overline \eta$, and $G_k\defeq \Gal (\overline k/k)$ the absolute Galois group of $k$.
Write
$$\Delta\defeq \pi_1(\overline {X},\overline \eta)\ \ \ \ \ \ \ \ \ \  \text {and}\ \ \ \ \ \ \ \ \ \ \Pi\defeq \pi_1(X, \eta).$$
Thus, we have a natural exact sequence 
$$ 1\to \Delta \to \Pi\to G_k\to 1.\tag 0.1$$

For an exact sequence $1\to H'\to H@>\pr>> G\to 1$ of profinite groups we will refer to a continuous homomorphism
$s:G\to H$ such that $\pr \circ s=\id_G$ as a group-theoretic {\it section}, or simply {\it section}, of the natural projection $\pr:H\twoheadrightarrow G$.

Suppose that the above exact sequence $(0.1)$ {\it splits}
(for example assume that $X(k)\neq \emptyset$). Let 
$$s:G_k\to \Pi$$ 
be a {\it section} of the above sequence. We view $\Delta$ as a (non-abelian) $G_k$-module
via the conjugation action of $s(G_k)$.

For the rest of $\S0$ we assume that $k$ is a {\it number field}, i.e., $k$ is a finite extension of $\Bbb Q$.
Let $v$ be a {\it prime} of $k$, $k_v$ the completion of $k$ at $v$, and $G_{k_v}\subset G_k$ a {\it decomposition group} associated to $v$. 
Thus, $G_{k_v}$ is only defined up to conjugation. We view $\Delta$ as a $G_{k_v}$-module
via the conjugation action of $s(G_{k_v})$. For each prime $v$ of $k$ we have a natural {\it restriction} map (of pointed non-abelian cohomology sets) 
$$\Res_v:H^1(G_k,\Delta)\to H^1(G_{k_v},\Delta),$$
and a natural map
$$\prod _{\all v} \Res _v: H^1(G_k,\Delta)\to \prod _{\all v} H^1(G_{k_v},\Delta),$$
where the product is over {\it all} primes $v$ of $k$. The {\it main  problem} we are concerned with in this paper is the following.

\definition {Question A} Is the map 
$$\prod _{\all v} \Res _v: H^1(G_k,\Delta)\to \prod _{\all v} H^1(G_{k_v},\Delta)$$
{\bf injective}?
\enddefinition

The above question is related to the {\it Grothendieck anabelian section conjecture} {\bf (GASC)}, which predicts that all splittings of the exact sequence
(0.1) arise from $k$-rational points of $k$ (See for example [Sa\"\i di], $\S0$, for a precise statement
of this conjecture). More precisely, there is a natural set-theoretic map 
$$X(k)\to H^1(G_k,\Delta),$$
and the {\bf (GASC)} predicts that
this map is {\it bijective}. The validity of the {\bf (GASC)} would then imply a positive answer to the above question. More precisely, the {\bf (GASC)} would imply 
the following (much stronger statement): the map 
$$H^1(G_k,\Delta)\to H^1(G_{k_v},\Delta)$$ 
is {\it injective if $v$ is a non-archimedean prime}. (Indeed, in this case  the map $X(k)\to H^1(G_k,\Delta)$ is bijective, 
the natural map  $X(k_v)\to H^1(G_{k_v},\Delta)$ is (well-known to be) injective,  
the image of the map 
$H^1(G_k,\Delta)\to H^1(G_{k_v},\Delta)$ would be contained in $X(k_v)$, and the natural map $X(k)\to X(k_v)$ is injective). 
Conversely, a positive answer to the above question
can be viewed as a kind of {\it (weak) positive evidence for the validity of the section conjecture}. 

Let $\Delta^{\sol}$ be the maximal {\bf prosolvable} quotient of $\Delta$, which is a {\it characteristic} quotient.
Then the above $G_k$ (resp. $G_{k_v}$)-module structure on $\Delta$ induce naturally a  $G_k$ (resp. $G_{k_v}$)-module structure on
$\Delta^{\sol}$. Let $S\subseteq \Primes _k$ be a {\it non-empty set of primes} of $k$, then we have as above a natural map

$$\prod _{v\in S} \Res _v^{\sol}: H^1(G_k,\Delta^{\sol})\to \prod _{v\in S} H^1(G_{k_v},\Delta^{\sol}).$$ 

Our main result in this paper is the following.
\proclaim {Theorem B} Assume $k$ is a {\bf number field}, and $S\subseteq \Primes _k$ is a set of primes of $k$ of {\bf density $1$}.
Then the map
$$\prod _{v\in S} \Res _v^{\sol}: H^1(G_k,\Delta^{\sol})\to \prod _{v\in S} H^1(G_{k_v},\Delta^{\sol})$$
is {\bf injective}.
\endproclaim

Our proof of Theorem B relies on a {\it devissage} argument and a careful analysis of the structure of the {\it geometrically prosolvable} arithmetic fundamental group of $X$,
and its splittings, which is established in $\S1$. The {\it abelian} analog of Theorem B is essentially a consequence of a Theorem of Serre (cf. Proposition 2.2). 
Theorem B is proved in $\S2$.

 \definition {Acknowledgment} I am very grateful to the referee for his/her valuable comments which helped improve an earlier version of this paper, especially for his/her comments
 concerning the proof of Proposition 2.2.
 \enddefinition

\subhead
\S1. Geometrically prosolvable arithmetic fundamental groups
\endsubhead
We use the notations in $\S0$. In particular,  $k$ is a field of characteristic $0$, and $X\to \Spec k$ a proper, smooth, and geometrically connected {\bf hyperbolic
curve} (i.e., $\genus (X)\ge 2$).  For a profinite group $H$ we denote by $\overline {[H,H]}$ the closed subgroup of $H$ which is {\it topologically generated by the commutator subgroup} of $H$.
Consider the {\it derived series} of $\Delta$
$$.....\subseteq \Delta (i+1)\subseteq \Delta (i)\subseteq......\subseteq \Delta (1)\subseteq \Delta(0)=\Delta, \tag 1.1$$
where 
$$\Delta(i+1)=\overline {[\Delta (i),\Delta(i)]},$$
for $i\ge 0$, is the $i+1$-th {\it derived subgroup}, which is a {\it characteristic} subgroup of $\Delta$. Write
$$\Delta _i\defeq \Delta/\Delta(i).$$
Thus, $\Delta _i$ is the {\bf $i$-th step prosolvable} quotient of $\Delta$, and $\Delta _1$ is the maximal {\bf abelian} quotient of $\Delta$. 
Note that there exists a natural exact sequence
$$ 1\to   \Delta^{i+1}\to \Delta _{i+1}\to \Delta _i\to 1 \tag 1.2$$
where $\Delta^{i+1}$ is the subgroup $\Delta (i)/\Delta (i+1)$ of $\Delta _{i+1}$. In particular, $\Delta ^{i+1}$ is {\it abelian}.

Write 
$$\Pi_i\defeq \Pi/\Delta (i),$$
which inserts in the following exact sequence 
$$1\to \Delta_i\to \Pi_i\to G_k\to 1.$$
We will refer to $\Pi_i$ as the {\bf geometrically $i$-th step prosolvable} fundamental group of $X$.
We have natural commutative diagrams
of exact sequences

$$
\CD
1  @>>>  \Delta   @>>> \Pi  @>>> G_k  @>>> 1  \\
@.   @VVV       @VVV       @VVV \\
1  @>>> \Delta _i  @>>>  \Pi_i   @>>>  G_k  @>>> 1 \\
\endCD
\tag 1.3
$$

and

$$
\CD
@.  1  @.  1  \\
@. @VVV    @VVV\\
@. \Delta ^{i+1}   @=  \Delta ^{i+1}  \\
@.  @VVV  @VVV\\
1   @>>>   \Delta _{i+1} @>>>  \Pi_{i+1}  @>>> G_k @>>> 1\\
@.       @VVV            @VVV             @VVV \\
1 @>>> \Delta _i   @>>> \Pi_i  @>>>  G_k @>>> 1\\
@.     @VVV       @VVV        @VVV \\
   @. 1 @. 1@. 1@. \\
\endCD
\tag 1.4
$$   
where the left and middle vertical maps in diagram (1.3) are natural surjections.

Write $\Delta ^{\sol}$ for the maximal {\bf prosolvable} quotient of $\Delta$, which is a {\it characteristic} quotient (cf. Lemma 1.1).
Consider the {\it pushout} diagram
$$
\CD
1 @>>> \Delta  @>>>  \Pi @>>> G_k  @>>> 1\\
@. @VVV      @VVV    @VVV  \\
1 @>>>  \Delta ^{\sol}   @>>> \Pi^{(\sol)}   @>>> G_k @>>> 1\\
\endCD
\tag 1.5
$$
We will refer to $\Pi^{(\sol)}$ as the {\bf geometrically prosolvable} fundamental group of $X$.

\proclaim {Lemma 1.1} We have natural identifications $\Delta ^{\sol}  \isom \underset{i\ge 1} \to{\varprojlim}\ \Delta_i$,
and $\Pi ^{(\sol)}\isom \underset{i\ge 1} \to{\varprojlim}\  \Pi_i$. 
In particular, $\Ker (\Delta \twoheadrightarrow \Delta ^{\sol})=\bigcap _{i\ge 1}\Delta (i)$,
and $\Delta ^{\sol}$ is a characteristic quotient of $\Delta$.
\endproclaim

\demo{Proof} Follows from the various definitions.
\qed
\enddemo

Let $i\ge 1$ be an integer. The profinite group $\Delta _{i}$ is {\it finitely generated} (as follows from the well-known
finite generation of $\Delta$ which projects onto $\Delta _{i}$).  
Let $\{\widehat \Delta _{i}[n]\}_{n\ge 1}$ be a {\it countable} system of {\bf characteristic open} subgroups of $\Delta_i$ such that 
$$\widehat \Delta_{i}[n+1]\subseteq \widehat \Delta _{i}[n],\ \ \ \ \widehat \Delta _{i}[1]\defeq \Delta_i,\ \ \ \ \text {and} \ \ \bigcap _{n\ge 1}\widehat \Delta_{i}[n]=\{1\}.$$
Write $\Delta _{i,n}\defeq \Delta_i/\widehat \Delta_i[n]$. Thus, $\Delta _{i,n}$ is a {\it finite characteristic} quotient of $\Delta _i$, 
which is an {\it i-th step solvable} group, and
we have a {\it pushout} diagram of exact sequences

$$
\CD
1@>>> \Delta_i @>>>   \Pi_i  @>>>  G_k @>>> 1\\
@.  @VVV    @VVV    @VVV \\
1@>>>  \Delta _{i,n}  @>>>  \Pi _{i,n} @>>>  G_k @>>> 1\\
\endCD
\tag 1.6
$$
which defines a ({\it geometrically finite}) quotient $\Pi_{i,n}$ of $\Pi_i$.

In the following discussion we {\it fix} an integer $i\ge 1$. Suppose that the exact sequence $1\to \Delta_i\to \Pi_i\to G_k\to 1$ {\bf splits}, and let
$$s_i:G_k\to \Pi_i$$
be a section of the natural projection $\Pi_i\twoheadrightarrow G_k$,
which induces a section
$$s_{i,n}:G_k\to \Pi_{i,n}$$ 
of the natural projection $\Pi_{i,n}\twoheadrightarrow G_k$ (cf. diagram (1.6)), for each $n\ge 1$.
Write 
$${\widehat \Pi _i}[n] \defeq {\widehat \Pi _i}[n][s_i] \defeq \widehat \Delta _i[n]. s_i (G_k).$$

Thus, ${\widehat \Pi_i}[n]\subseteq \Pi_i$ is an open subgroup which contains the image 
$s_i(G_k)$ of $s_i$. Write 
$${\Pi_i}[n]\defeq {\Pi_i}[n][s_i]$$ 
for the {\it inverse image} of $\widehat \Pi _i[n]$
in $\Pi$ (cf. diagram (1.3)). Thus, $\Pi_i[n]\subseteq \Pi$ is an {\it open} subgroup which corresponds
to an \'etale cover 
$$X_{i,n}\to X$$
where $X_{i,n}$ is a geometrically irreducible $k$-curve (since $\Pi_i[n]$ maps onto $G_k$ via the natural projection $\Pi\twoheadrightarrow G_k$, by the very 
definition of $\Pi_i[n]$). 
Note that the \'etale cover $\overline X_{i,n}\defeq X_{i,n}\times _k\overline k\to \overline X$ is Galois with
Galois group $\Delta _{i,n}$, and we have a commutative diagram of \'etale covers

$$
\CD
\overline X_{i,n}  @>>>   \overline X \\
@VVV    @VVV \\
X_{i,n}  @>>>  X\\
\endCD
$$
where $\overline X_{i,n}\to X$ is Galois with Galois group $\Pi_{i,n}$, and $\overline X_{i,n}\to X_{i,n}$ is Galois with Galois group $s_{i,n}(G_k)$. 
Moreover, we have a commutative diagram of exact sequences
$$
\CD
@.   1@.    1@.\\
@. @VVV  @VVV\\
1 @>>>  \Delta _i[n]= \pi_1(\overline X_{i,n},\overline \eta_{i,n})@>>> \Pi_i[n]=\pi_1(X_{i,n},\eta_{i,n}) @>>>  G_k @>>> 1\\
@.   @VVV    @VVV   @VVV \\
1 @>>> \Delta=\pi_1(\overline {X},\overline \eta) @>>>   \Pi=\pi_1(X,\eta) @>>>  G_k  @>>> 1 \\
\endCD
\tag 1.7
$$
where $\Delta_i[n]$ is the inverse image of $\widehat \Delta_i[n] $
in $\Delta $, and the equalities $\Delta _i[n]= \pi_1(\overline X_{i,n},\overline \eta_{i,n})$, $\Pi_i[n]=\pi_1(X_{i,n},\eta_{i,n})$,
are natural identifications. Note that since the \'etale cover $X_{i,n}\to X$ (resp. $\overline X_{i,n}\to \overline X$) 
is defined via an open subgroup of $\pi_1(X,\eta)$ (resp. $\pi_1(\overline {X},\overline \eta)$) it is a pointed \'etale cover, and   $X_{i,n}$ (resp. $\overline X_{i,n}$) is naturally endowed with a geometric point $\eta_{i,n}$ (resp. $\overline \eta_{i,n}$) above $\eta$ (resp. $\overline \eta$).
Note also that $\Pi_i[n+1]\subseteq \Pi_i[n]$, and $\Delta _i[n+1]\subseteq \Delta_i[n]$, 
as follows from the various definitions.

For each integer $n\ge 1$, write $\Delta _i[n]^{\ab}$ for the maximal {\it abelian} quotient of $\Delta _i[n]$, which is a {\it characteristic} quotient,
and consider the {\it pushout} diagram
$$
\CD
1 @>>>  \Delta _i[n] @>>> \Pi_i[n] @>>>  G_k @>>> 1\\
@.  @VVV    @VVV   @VVV \\
1 @>>>  \Delta _i[n]^{\ab} @>>> \Pi_i[n]^{(\ab)} @>>>  G_k @>>> 1\\
\endCD
\tag 1.8
$$
Thus, $\Pi_i[n]^{(\ab)}$ is the {\it geometrically abelian} fundamental group of $X_{i,n}$. Further, consider the following commutative diagram

$$
\CD
1 @>>>  \Delta^{i+1}   @>>>  \Cal H_i\defeq \Cal H_{i}[s_i] @>>>  G_k  @>>> 1\\
@.    @VVV     @VVV     @V{s_i}VV \\
1  @>>>   \Delta ^{i+1}  @>>>   \Pi _{i+1}   @>>>  \Pi_i @>>> 1\\
\endCD
\tag 1.9
$$
where the lower exact sequence is the sequence in diagram (1.4), and the right square is {\bf cartesian}. Thus, (the group extension) $\Cal H_i$
is the {\it pullback} of (the group extension) $\Pi_{i+1}$ via the section $s_i$.

\proclaim {Lemma 1.2} We have natural identifications $\Delta ^{i+1}  \isom \underset{n\ge 1} \to{\varprojlim}\ \Delta_i[n] ^{\ab}$,
and
$\Cal H_i\isom \underset{n\ge 1} \to{\varprojlim}\  \Pi_i[n]^{(\ab)}$.
\endproclaim

\demo{Proof} Follows from the various definitions.
\qed
\enddemo

Next, assume that the above setion $s_i:G_k\to \Pi_i$ can be {\it lifted} to a section $s_{i+1}:G_k\to \Pi_{i+1}$
of the natural projection $\Pi_{i+1}\twoheadrightarrow G_k$, i.e., there exists such a section $s_{i+1}$ and a commutative diagram
$$
\CD
G_k  @>{s_{i+1}}>> \Pi_{i+1}\\
@V{\id}VV     @VVV\\
G_k  @>{s_i}>> \Pi _i\\
\endCD
$$
where the right vertical map is as in diagram (1.4). Recall the exact sequence (1.2)
$$1\to \Delta ^{i+1}\to \Delta _{i+1}\to \Delta _i\to 1.$$ 
We view $\Delta _{i+1}$
(hence also $\Delta _i$, and $\Delta ^{i+1}$) as a $G_k$-module via the action of $s_{i+1}(G_k)$ 
by conjugation. Thus, this $G_k$-module structure on $\Delta _i$ is the one induced by the section $s_i$.
The above sequence is an exact sequence of $G_k$-modules.
\proclaim {Lemma 1.3} Assume that $k$ is a {\bf number field}, or {\bf a $p$-adic local field} 
(i.e., a finite extension of $\Bbb Q_p$). Then $H^0(G_k,\Delta ^{i+1})=\{0\}$, $\forall i\ge 0$.
\endproclaim

\demo{Proof} 
The group $H^0(G_k,\Delta ^{i+1})$ is naturally identified with $\underset{n\ge 1} \to{\varprojlim}\ H^0(G_k,\Delta_i[n] ^{\ab})$
(cf. Lemma 1.2), where the $G_k$-module $\Delta_i[n] ^{\ab}$ is isomorphic to the Tate module of the jacobian of the curve $X_{i,n}$.
Furthermore, $H^0(G_k,,\Delta_i[n] ^{\ab})=\{0\}$, as follows immediately from the well-known fact that the torsion subgroup of the group of $k$-rational points $A(k)$  of an abelian variety $A$ over $k$ is finite, which is a consequence  of the Mordell-Weil Theorem in the case $k$ is a number field (cf. [Serre], 4.3), and [Mattuck], Theorem 7, p.114, in the case where $k$ is a $p$-adic local field.
\qed
\enddemo

\proclaim{Lemma 1.4} Assume that $k$ is a {\bf number field}, or {\bf a $p$-adic local field} (i.e., a finite extension of $\Bbb Q_p$). Then, $H^0(G_k,\Delta_i)=\{1\}$, $\forall i\ge 0$.
\endproclaim

\demo{Proof} By induction on $i$, using Lemma 1.3, and the fact that we have an exact sequence of groups
$1\to H^0(G_k,\Delta ^{i+1})\to H^0(G_k,\Delta _{i+1})\to H^0(G_k,\Delta _i)$.
\qed
\enddemo

\proclaim {Lemma 1.5} Assume that $k$ is {\bf number field}, or {\bf a $p$-adic local field} (i.e., a finite extension of $\Bbb Q_p$). Then the natural map
$H^1(G_k,\Delta ^{i+1})\to H^1(G_k,\Delta _{i+1})$ of pointed sets is {\bf injective}.
\endproclaim

\demo{Proof} There exists an exact sequence of pointed sets 
$$H^0(G_k,\Delta _i)\to H^1(G_k,\Delta ^{i+1})\to H^1(G_k,\Delta _{i+1})$$
(cf. [Serre1], I, $\S5$, 5.5, Proposition 38). The proof then follows from [Serre1], I, $\S5$, 5.5, Proposition 39 (ii), and
the fact that $H^0(G_k,\Delta_i)$ is trivial (cf. Lemma 1.4).
\qed
\enddemo

\subhead
\S2. Proof of Theorem B
\endsubhead
This section is devoted to the proof of Theorem B (cf. $\S0$).
We use the same notations as in Theorem B, $\S0$, and $\S1$. We assume further that the set $S\subset \Primes_k$ 
{\it contains no real places}.

Recall the commutative diagram of exact sequences of profinite groups 
$$
\CD
1@>>> \Delta \defeq \pi_1(\overline {X},\overline \eta) @>>> \Pi\defeq \pi_1(X, \eta) @>>> G_k @>>> 1\\
@.  @VVV   @VVV   @VVV\\
1@>>> \Delta ^{\sol}@>>>   \Pi^{(\sol)}@>>> G_k@>>> 1\\
\endCD
$$
(cf. diagram (1.5)). And the natural map (cf $\S0$)
$$\prod _{v\in S} \Res _v^{\sol}: H^1(G_k,\Delta^{\sol})\to \prod _{v\in S} H^1(G_{k_v},\Delta ^{\sol}).$$ 
We will show this map is {\bf injective}. Recall the definition of the {\it $i$-th step prosolvable characteristic} quotient $\Delta _i$ of $\Delta$ 
(cf. the discussion before the exact sequence (1.2)).

\proclaim {Proposition 2.1}
The natural map 
$$\prod _{v\in S} \Res _{v}^i: H^1(G_k,\Delta_i)\to \prod _{v\in S} H^1(G_{k_v},\Delta _i)$$ 
is {\bf injective}.
\endproclaim

We will prove Proposition 2.1 by an induction argument on $i\ge 1$. The case $i=1$ follows from the following.

\proclaim {Proposition 2.2}
Let $A\to \Spec k$ be an {\bf abelian variety over a number field $k$}, and $TA$
its Tate module. Let $S\subset \Primes _k$ be a set of primes of $k$ of {\bf density $1$}.
Then the natural homomorphism
$$\prod _{v\in S} \Res _{v}: H^1(G_k,TA)\to \prod _{v\in S} H^1(G_{k_v},TA)$$ 
is {\bf injective}.
\endproclaim

\demo{Proof} This follows from a Theorem of Serre and is essentially proven in [Milne], Chapter I,  Proposition 6.22, page 110. 
The proof in loc. cit. however contains an inaccuracy  
related to the use of Chebotarev's density Theorem which we fix in the following.
It suffices to show, given a prime number $l$, that the homomorphism 
$$\prod _{v\in S} \Res _{v}: H^1(G_k,T_lA)\to \prod _{v\in S} H^1(G_{k_v},T_lA)$$ 
is injective,
where $T_lA$ is the $l$-adic Tate module.
For each integer $n\ge 1$, let $k_n\defeq k(A[l^n])$ be the field obtained by adjoining to $k$ the coordinates of the $l^n$-torsion points of $A$,
and $k_{\infty}\defeq \bigcup_{n\ge 1} k_n$. Let $S_{n}\subseteq \Primes_{k_{n}}$ be the pre-image of $S$ in $\Primes_{k_n}$ via the natural map 
$\Primes_{k_n}\to \Primes_k$, which is a set of primes of $k_{n}$ of density $1$.
Let $G_{k_n}\defeq \Gal (\bar k/k_n)$. Then the natural restriction homomorphism
$$\prod _{w\in S_n} \Res _{w}: H^1(G_{k_n},A[l^n])\to \prod _{w\in S_n} H^1(G_{k_n,w},A[l^n]),$$ 
where $G_{k_n,w}$ is a decomposition subgroup of $G_{k_n}$ associated to $w$, and $A[l^n]$ is the Galois module of $l^n$-torsion points of $A(\bar k)$, is injective. 
Indeed, the groups $H^1(G_{k_n}, A[l^n])$, and $H^1(G_{k_n,w},A[l^n])$, are groups of continuous homomorphisms $\Hom(G_{k_n},A[l^n])$, and $\Hom(G_{k_n,w},A[l^n])$, respectively, and the map 
$$\prod _{w\in S_n} \Res _{w}: \Hom (G_{k_n},A[l^n])\to \prod _{w\in S_n} \Hom (G_{k_n,w},A[l^n])$$ 
is injective by Chebotarev's density Theorem (cf. [Neukirch-Schmidt-Winberg], Chapter IX, (9.1.3) Theorem (i)). 
Therefore 
$$L_n\defeq \Ker \lgroup  H^1(G_{k},A[l^n])@>\prod _{v\in S} \Res _{v} >> \prod _{v\in S} H^1(G_{k_v},A[l^n] \rgroup$$ 
is contained in the image of the inflation map 
$H^1(\Gal (k_n/k),A[l^n]) @>\inf>> H^1(G_k,A[l^n])$ which factorizes as 
$H^1(\Gal (k_n/k),A[l^n]) @>{\inf}>> H^1(\Gal (k_{\infty}/k),A[l^n]) @>{\inf}>> H^1(G_k,A[l^n])$.
From this we deduce that 
$$L\defeq \Ker \lgroup H^1(G_k,T_lA)@>\prod _{v\in S} \Res _{v}>>
\prod _{v\in S} H^1(G_{k_v},T_lA)\rgroup$$ 
is naturally identified with $\underset{n\ge 1} \to{\varprojlim} L_n$
(cf. [Neukirch-Schmidt-Winberg], Chapter II, (2.3.5)Corollary, and the fact that ${\varprojlim}$ is left exact), and $L$
is contained in the image of the inflation map  $\inf:H^1(\Gal (k_{\infty}/k),T_lA)\to H^1(G_k,T_lA)$. 

Now $H^1(\Gal (k_{\infty}/k),T_lA)$ is a finite group by a Theorem of Serre (cf. [Milne], Chapter I, Proposition 6.19, page 109). Thus, 
$L$ is a torsion group. Recall the Kummer exact sequence
$0\to A(k)^l \to H^1(G,T_lA)\to T_lH^1(G,A)\to 0$, where $T_lH^1(G,A)$ is the $l$-adic Tate module of $H^1(G,A)$, and $A(k)^l$ is the $l$-adic completion of $A(k)$.
The group $L$ is contained in $A(k)^l$ since $T_lH^1(G_k,A)$ is torsion free, and is further contained in the torsion subgroup $A(k)^{l, \tor}$ 
 of $A(k)^{l}$. The map $A(k)^{l,\tor}\to \prod_{v\in S} A(k_v)^{l,\tor}$, where $A(k_v)^{l,\tor}$ is the torsion subgroup of the $l$-adic completion $A(k_v)^l$ of $A(k_v)$,
is injective. Hence $L$ is trivial.
\qed
\enddemo

\demo {Proof of Proposition 2.1}
Fix an integer $i\ge 1$.  Consider the following commutative diagram of maps of pointed cohomology sets
$$
\CD
1 @>>> H^1(G_k,\Delta ^{i+1}) @>>> H^1(G_k,\Delta _{i+1}) @>>> H^1(G_k,\Delta _{i})  \\
@.   @VVV   @V{\prod _{v\in S} \Res _{v}^{i+1}}VV   @V{\prod _{v\in S} \Res _{v}^i}  VV \\
1 @>>> \prod _{v\in S} H^1(G_{k_v},\Delta ^{i+1}) @>>> \prod _{v\in S} H^1(G_{k_v},\Delta _{i+1})@>>> \prod _{v\in S} H^1(G_{k_v},\Delta _{i})\\
\endCD
$$
where the horizontal sequences are exact (cf. Lemma 1.5) and the vertical maps are the natural restriction maps.
We assume by induction hypothesis that the right vertical map
$\prod _{v\in S} \Res _{v}^i: H^1(G_k,\Delta_i)\to \prod _{v\in S} H^1(G_{k_v},\Delta _i)$
is {\bf injective}. We will show that the middle vertical map
$\prod _{v\in S} \Res _{v}^{i+1}: H^1(G_k,\Delta_{i+1})\to \prod _{v\in S} H^1(G_{k_v},\Delta _{i+1})$ 
is {\bf injective}.
Let 
$$[\rho], [\tau]\in H^1(G_k,\Delta_{i+1})$$
be two cohomology classes such that
$$\prod _{v\in S} \Res _{v}^{i+1}([\tau])=\prod _{v\in S} \Res _{v}^{i+1}([\rho]).$$
Write $s_i:G_k\to \Pi_{i}$ for the section of the projection $\Pi_{i}\twoheadrightarrow G_k$ induced by the section $s:G_k\to \Pi$, for $i\ge 1$.
We will show $[\tau]=[\rho]$. We can (without loss of generality), and will, assume that $[\tau]=[s_{i+1}]=1$ is the distinguished element of $H^1(G_k,\Delta _{i+1})$.
The classes $[\rho]$, and $[\tau]$, map to the classes $[\rho_1]$, and $[\tau_1]=[s_i]=1$ in $H^1(G_k,\Delta _i)$, respectively.
In particular, we have the equality
$$\prod _{v\in S} \Res _{v}^i([\tau_1])=\prod _{v\in S} \Res _{v}^i([\rho_1]),$$
hence $[\tau_1]=[\rho_1]=1$, since $\prod _{v\in S} \Res _{v}^i$ is injective by the induction hypothesis. 
Thus, there exist classes $[\tilde \tau]=1$, and $[\tilde \rho]$, in $H^1(G_k,\Delta ^{i+1})$ which map to $[\tau]$, and $[\rho]$, respectively
in $H^1(G_k,\Delta _{i+1})$ (cf. exactness of horizontal sequences in the above diagram).

Next, and in order to show that $[\rho]=[\tau]$ in $H^1(G_k,\Delta_{i+1})$,
it suffices to show $[\tilde \rho]=[\tilde \tau]$ in $H^1(G_k,\Delta^{i+1})$.
Note that the assumption $\Res_v^{i+1}([\tau])= \Res_v^{i+1}([\rho])$
implies that  $\Res_v^{i+1} ([\tilde \tau])= \Res_v^{i+1}([\tilde \rho])$ in $H^1(G_{k_v},\Delta ^{i+1})$, for each place $v\in S$,
as follows from the injectivity of the maps $H^1(G_{k_v},\Delta^{i+1})\to H^1(G_{k_v},\Delta_{i+1})$
(cf. Lemma 1.5). 
We have a commutative diagram of group homomorphisms
$$
\CD
H^1(G_k,\Delta^{i+1})  @>>> \underset{n\ge 1} \to{\varprojlim} H^1(G_k, \Delta _i[n]^{\ab})\\
@V{\prod _{v\in S} \Res _{v}^{i+1}}VV                      @VVV\\
\prod _{v\in S}H^1(G_{k_v},\Delta ^{i+1})        @>>>     \prod _{v\in S}  \lgroup \underset{n\ge 1} \to{\varprojlim}  H^1(G_{k_v},\Delta _i[n]^{\ab})\rgroup\\
\endCD
$$
where the horizontal maps are induced from the natural identification $\Delta ^{i+1}  \isom \underset{n\ge 1} \to{\varprojlim}\ \Delta_i[n] ^{\ab}$ (cf. Lemma 1.2). 
This identification induces natural isomorphisms
$H^1(G_k,\Delta^{i+1})\isom \underset{n\ge 1} \to{\varprojlim} H^1(G_k, \Delta _i[n]^{\ab})$,
and $H^1(G_{k_v},\Delta^{i+1})\isom \underset{n\ge 1} \to{\varprojlim} H^1(G_{k_v}, \Delta _i[n]^{\ab})$,
for each $v\in S$, where the transition morphisms in the projective limit are the natural homomorphisms induced by the natural
$G_k$ (resp. $G_{k_v}$)-homomorphisms $\Delta _{i}[n+1]^{\ab}\to \Delta _{i}[n]^{\ab}$ (cf. [Neukirch-Schmidt-Winberg], Chapter II, (2.3.5)Corollary). 
The right vertical map in the above diagram is injective by Lemma 2.2 
(recall the $G_k$-module $\Delta_i[n] ^{\ab}$ is isomorphic to the Tate module of the jacobian of the curve $X_{i,n}$ defined in $\S1$). 
Thus, the left vertical map
$H^1(G_k,\Delta ^{i+1})\to \prod _{v\in S}H^1(G_{k_v},\Delta ^{i+1})$
is {\bf injective}. In particular, $[\tilde \rho]=[\tilde \tau]$, and $[\rho]=[\tau]$.
This finishes the proof of Proposition 2.1.
\qed
\enddemo

Next, let 
$$[\alpha], [\beta]\in H^1(G_k,\Delta^{\sol})$$
be two cohomology classes such that
$$\prod _{v\in S} \Res _{v}^{\sol}([\alpha])=\prod _{v\in S} \Res _{v}^{\sol}([\beta]).$$
We will show $[\alpha]=[\beta]$. The above cohomology classes $[\alpha], [\beta]\in H^1(G_k,\Delta^{\sol})$
give rise to cohomology classes 
$$([\alpha_i])_{i\ge 1}, ([\beta_i])_{i\ge 1}\in \underset {i\ge1} \to{\varprojlim}\ H^1(G_k,\Delta_i),$$
(cf. Lemma 1.1). 
For $i\ge 1$, the following holds (by assumption) 
$$\prod _{v\in S} \Res _{v}^{i}([\alpha_i])=\prod _{v\in S} \Res _{v}^{i}([\beta_i]).$$
Thus,
$$[\alpha_i]=[\beta_i],\ \ \ \ \ \ \ \forall i\ge 1,$$
by Proposition 2.1. Let $\alpha:G_k\to \Pi^{(\sol)}$, and $\beta:G_k\to \Pi^{(\sol)}$, be sections corresponding 
to $[\alpha]$, and $[\beta]$, respectively.
We will show that the sections $\alpha$ and $\beta$ are conjugated by some element of $\Delta ^{\sol}$. Let
$\alpha_i:G_k\to \Pi_i$, and $\beta_i:G_k\to \Pi_i$, be the sections induced by $\alpha$, and $\beta$, respectively.

\proclaim {Lemma 2.3} For each $i\ge 1$, there exists $\sigma_i\in \Delta _i$ such that $\sigma_i\alpha_i\sigma_i^{-1}=\beta_i$ (equality of sections of the projection $\Pi_i\twoheadrightarrow G_k$), 
and $\sigma_{i}$ maps to $\sigma_j$ via the natural homomorphism $\Delta_{i}\twoheadrightarrow \Delta _j$, for $j\le i$. Thus, $(\sigma_i)_{i\ge 1}\in {\underset{i\ge 1} 
\to{\varprojlim}}\  \Delta _i=\Delta ^{\sol}$ (cf. Lemma 1.1).
\endproclaim

\demo{Proof} For each $i\ge 1$, we have $[\alpha_i]=[\beta_i]$ in $H^1(G_k,\Delta_i)$ by assumption. Thus,    
there exists $\sigma_i\in \Delta _i$ such that $\sigma_i\alpha_i\sigma_i^{-1}=\beta_i$ (equality of sections of the projection $\Pi_i\twoheadrightarrow G_k$).
The set $S_i\defeq \{\tilde \sigma_i\in \Delta_i\ \ \vert\ \ \tilde \sigma_i\alpha_i\tilde \sigma_i^{-1}=\beta_i\}$ is non-empty and equals
$\sigma_i.N_{\Pi_i}(\alpha_i(G_k))\cap \Delta_i$, where $N_{\Pi_i}(\alpha_i(G_k))$ is the normaliser of $\alpha_i(G_k)$ in $\Pi_i$. 
In particular, $S_i$ is compact. The projective limit  $\underset {i\ge1} \to{\varprojlim}\ S_i$ is non-empty since the $\{S_i\}_{i\ge 1}$ 
are compact and non-empty (cf. [Ribes-Zalesskii], Chap. 1, Proposition 1.1.4). An element $\sigma=(\sigma_i)_{i\ge 1}\in \underset {i\ge1} \to{\varprojlim}\ S_i$ is as required in Lemma 2.3. 
\qed
\enddemo

Let $\sigma =(\sigma_i)_{i\ge 1}\in \Delta ^{\sol}$ be as in Lemma 2.3. The set $\Sec (G_k,\Pi ^{\sol})$  of  
sections of the projection $\Pi^{\sol}\twoheadrightarrow G_k$
is naturally in bijection with the projective limit $\underset {i\ge1} \to{\varprojlim}\  \Sec(G_k,\Pi_i)$ of the sets $\Sec(G_k,\Pi_i)$ of sections of the projection 
$\Pi_i\twoheadrightarrow G_k$ (cf. Lemma 1.1). Then it follows from Lemma 2.3 that $\sigma \alpha \sigma ^{-1}=\beta$ as elements of $\Sec (G_k,\Pi ^{\sol})$,
hence $[\sigma \alpha \sigma ^{-1}]=[\alpha]=[\beta]$ in $H^1(G_k,\Delta ^{\sol})$ as desired.

This finishes the Proof of Theorem B.
\qed

$$\text{References.}$$

\noindent
[Grothendieck] Grothendieck, A., Rev\^etements \'etales et groupe fondamental, Lecture 
Notes in Math. 224, Springer, Heidelberg, 1971.

\noindent
[Mattuck] Mattuck, A., Abelian varieties over $p$-adic ground fields, Annals of Mathematics, Second Series, Vol. 62, No. 1, 92-119.

\noindent
[Milne] Milne, J.S., Arithmetic duality theorems, Perspectives in Mathematics, Volume 1, J. Coates and  S. Helgason, editors, 1986.

\noindent
[Neukirch-Schmidt-Winberg] Neukirch, J., Schmidt, A., Winberg, K., Cohomology of number fields, first edition,
Springer, Grundlehren der mathematischen Wissenschaften Bd. 323, 2000.

\noindent
[Ribes-Zalesskii] Ribes, L., Zalesskii, P., Profinite Groups, Ergebnisse der Mathemtik und ihrer Grenzgebiete, 
Folge 3, Volume 40, Springer, 2000.

\noindent
[Sa\"\i di] Sa\"\i di, M., The cuspidalisation of sections of arithmetic fundamental groups, Advances in Mathematics 230
(2012) 1931-1954.

\noindent
[Serre] Serre, J-P., Lectures on the Mordell-Weil Theorem, Translated and edited by Martin Brown from notes by Michel Waldschmidt, $2^{\nd}$ edition, 
Aspect of Mathematics, 1990.

\noindent
[Serre1] Serre, J-P., Galois cohomology, Springer-Verlag Berlin Heidelberg, 1997.

\bigskip

\noindent
Mohamed Sa\"\i di

\noindent
College of Engineering, Mathematics, and Physical Sciences

\noindent
University of Exeter

\noindent
Harrison Building

\noindent
North Park Road

\noindent
EXETER EX4 4QF 

\noindent
United Kingdom

\noindent
M.Saidi\@exeter.ac.uk

\end
\enddocument